\documentclass[12pt]{article}
\usepackage{amsmath,amsthm,amssymb,amscd}
\usepackage[dvips]{graphicx}
\usepackage{here}
\usepackage{latexsym}
\title{On the detectability of different forms of
interaction in regression models}
\author{Juxin LIU$^1$  and Paul GUSTAFSON$^2$\\
\ \\
$^1$ Department of Medicine, University of Alberta\\
juxin@ualberta.ca\\
$^2$ Department of Statistics, University of British Columbia\\
 gustaf@stat.ubc.ca}
\date{\today}
\begin{document}

\def\bX{\mathbf{X}}
\def\bx{\mathbf{x}}
\def\bz{\mathbf{z}}
\def\bbX{\mathbb{X}}
\def\bbY{\mathbb{Y}}
\def\bY{\mathbf{Y}}
\def\bT{\mathbf{T}}
\def\bS{\mathbf{S}}
\def\bbS{\mathbb{S}}
\def\dd{{\rm{d}\,}}
\def\EE{{\rm{E}\,}}
\def\VV{{\rm{Var}\,}}
\def\Cov{{\rm{Cov}}}
\def\tr{{\rm{tr}}}
\def\bW{\mathbf{W}}
\def\bI{\mathbf{I}}
\def\bJ{\mathbf{J}}
\def\bu{\mathbf{u}}
\def\bv{\mathbf{v}}

\def\Beta{\boldsymbol{\beta}}
\def\Alpha{\boldsymbol{\alpha}}
\def\bTheta{\boldsymbol{\theta}}
\def\Epsilon{\boldsymbol{\epsilon}}
\def\bomega{\boldsymbol{\omega}}
\def\bgamma{\boldsymbol{\gamma}}
\def\bPhi{\boldsymbol{\Phi}}
\def\bsigma{\boldsymbol{\sigma}}
\def\bfeta{\boldsymbol{\eta}}
\def\bphi{\boldsymbol{\phi}}
\def\blamb{\boldsymbol{\lambda}}
\def\bzeta{\boldsymbol{\zeta}}

\newlength{\single}
\setlength{\single}{1.0\baselineskip}
\newlength{\double}
\setlength{\double}{1.4\baselineskip}

\baselineskip\double

\maketitle
\begin{abstract}
We derive an
asymptotic power function
for a likelihood-based test for interaction in a regression model,
with possibly misspecified alternative
distribution.
This allows a general investigation of types of interactions which are poorly or well
detected via data.
Principally we contrast pairwise-interaction models
with `diffuse interaction models' as introduced
in Gustafson, Kazi, and Levy (2005).
\end{abstract}

\section{Introduction}
There has been much discussion about how to define and measure
interaction. The interaction of two or more covariates can be
measured as the difference between the joint effect of covariates
and the sum of their independent effects, or, in other words, the
departure away from an additive model. In this paper, we focus on
the power of model-based tests for the presence of interaction,
under misspecified models. That is, with one kind of interaction
model truly generating the data, another kind of interaction model
is applied for estimation and testing purposes. Our rationale for
this is a suspicion that sometimes a wrong but parsimonious model
for interactions may lead to better power for detecting departures
from additivity than a complex model for interactions, even if the
complex model is correct. To make comparisons tractable, we derive
the asymptotic power of the test statistic under a sequence of
Pitman-type alternatives, which are getting closer to the null
(additive) model as the sample size increases.

In Section 2, under a general framework, we give the asymptotic
power function based on a misspecified model.  The corresponding
function for a correct model arises as a special case.  In Section
3, we apply the general results to a particular comparison between
a pairwise interaction model (PIM) and a diffuse interaction model
(DIM).  The latter was proposed by Gustafson, Kazi and Levy (2005)
as a parsimonious model for interactions appropriate for
reflecting a general synergism or antagonism in how covariates
interact, without identifying particular pairs of variables
responsible.  We find that when the DIM is correct, the DIM-based
test for interaction is more powerful than the PIM-based test, at
least in all the specific scenarios we have considered. When, the
PIM is correct, however, the comparison is mixed.  Depending on
the specific nature of the pairwise interactions, in terms of
directions and relative magnitudes of coefficients, either the
DIM-based test or the PIM-based test may be more powerful.

\section{General framework}

In this section, we give a general result about the
asymptotic power function of a Wald (quadratic form) test
for the presence of interactions, in the context of a misspecified
model for the alternative distribution.
In fact, the mathematical formalism is more general, in terms
of describing an arbitrary testing scenario with model misspecification.
Let $\mathcal{F}=\{f(y|\;\bx,\bTheta):\bTheta\in\Theta\}$
and $\mathcal{G}=\{g(y|\;\bx,\bomega), \bomega\in\Omega\}$ denote
two different parametric families of densities for modelling
$(Y|X_{1},\ldots,X_{p})$, with
$p_1=\mbox{dim}(\bTheta)$
and
$p_2=\mbox{dim}(\bomega)$.
We consider fitting model $\mathcal{F}$
to a sample of size $n$,
and testing
the null hypothesis $C\bTheta = \bzeta_{0}$ against a non-directional alternative,
where $C$ is an $r \times p_{1}$ matrix of full row rank.
Conversely,
the true data-generating mechanism is taken to be a member of $\mathcal{G}$.
To form a sequence of Pitman-type alternatives (Le Cam 1960),
the specific member of $\mathcal{G}$ generating the data is taken to
be
$\bomega_n=\bomega_0+n^{-1/2}\Delta\bfeta$,
where $\Delta$ is a scalar and $\bfeta$ is a vector of unit-length.
It is assumed that
$g(y| \bx, \bomega_{0}) \equiv f(y| \bx, \bTheta_{0})$ for
some $\bTheta_{0} \in \Theta$ with $C\bTheta_{0}=\bzeta_{0}$.
Thus the extent of model misspecification and the extent of deviation from the null both diminish with $n$.

For the two parametric families, let
$$
s_F(\bTheta,Y,\bX)=\partial{[\log\{f(Y|\;\bX,\bTheta)\}]}/\partial{\bTheta}
$$
 and
 $$
 s_G(\bomega,Y,\bX)=\partial{[\log\{g(Y|\;\bX,\bomega)\}]}/\partial{\bomega}
 $$
be the respective score vectors, and
$$
I_F(\bTheta)=\EE_{\bTheta}\{s_F(\bTheta,Y,\bX)s^T_F
(\bTheta,Y,\bX)\}
$$
and
$$
I_G(\bomega)=\EE_{\bomega}\{s_G(\bomega,Y,\bX)s^T_G
(\bomega,Y,\bX)\}
$$
be the respective Fisher information matrices.
Note that here, and in what follows, expectations are with respect
to the same fixed distribution of $\bX$, and the distribution of
$(Y|\bX)$ based on a member of $\cal F$ or $\cal G$, as indicated by
a subscript.
Let $\widehat\bTheta_n$ be the maximum likelihood estimator based on
fitting $\cal{F}$ to
$n$ observations arising from
$g(y|\;\bx,\bomega_n)$.
Then
$$W=n
\big(
  C\widehat\bTheta_n - \bzeta_{0}
\big)^T \left\{CI^{-1}_F(\widehat\bTheta_n)C^T \right\}^{-1} \big(
  C\widehat\bTheta_n - \bzeta_{0}
\big)
$$
is a Wald test
statistic,
which would be asymptotically distributed as $\chi_{r}^{2}$
with $r$=rank$(C)$ degrees of freedom if the data were generated under
the null (i.e., via some member of $\cal F$ with
$C\bTheta = \bzeta_{0}$).
However, with data generated as $g(y|\;\bx,\bomega_n)$, we have
\begin{eqnarray} \nonumber
n^{1/2}\left(C\widehat\bTheta - \bzeta_{0}\right)
& = &
n^{1/2}C(\widehat\bTheta_n-\bTheta_0) \\
& = & \label{truedist}
n^{1/2}C\left(\widehat\bTheta_n-\bTheta_*(\bomega_n)\right)
+n^{1/2}C\left(\bTheta_*(\bomega_n)-\bTheta_0\right),
\end{eqnarray}
where
$\bTheta_*(\bomega)$ is the parameter vector which minimizes
the Kullback-Leibler information criterion, that is
\begin{eqnarray} \label{kl}
\bTheta_*(\bomega)=\mbox{argmin}_{\bTheta}
E_{\omega}
\left\{\log
\frac{g(Y|\;\bX,\bomega)}{f(Y|\;\bX,\bTheta)}
\right\}.
\end{eqnarray}
Note that the fact $g(\cdot|\;\bomega_0)=f(\cdot|\;\bTheta_0)$
yields that $\bTheta_*(\bomega_0)=\bTheta_0$.

By White (1982) we know that the first item on the right side of
($\ref{truedist}$) is asymptotically normal with mean $\mathbf{0}$
and covariance matrix $CI^{-1}_F(\bTheta_0)C^T $. So we only need
to work on the second item. By (\ref{kl}) we know that
$\bTheta_*(\bomega)$ satisfies
\begin{eqnarray}\label{deriv}
E_{\omega}\left\{ s_F(\bTheta(\bomega),Y,\bX) \right\}
& = & 0.
\end{eqnarray}
Based on Gustafson (2001), implicit differentiation of
($\ref{deriv}$) gives
$$
\EE_{\bomega}[
s^{\prime}_F(\bTheta_*(\bomega);Y,\bX)]\frac{\partial\bTheta_*}{\partial\bomega}
+\EE_{\bomega}[s_F(\bTheta_*(\bomega),Y,\bX)s^{T}_G(\bomega;Y,\bX)]=0.
$$
Evaluated at $\bomega=\bomega_0$, the above equality yields
\begin{equation}\label{deriv1}
\frac{\partial\bTheta_*}{\partial\bomega}\Bigm
|_{\bomega=\bomega_0}=I^{-1}_F(\bTheta_0)\EE_{\bTheta_0}
\{s_F(\bTheta_0;Y,\bX)s^{T}_G(\bomega_0;Y,\bX)\},
\end{equation}
which is derived by the fact that
$\bTheta_*(\bomega_0)=\bTheta_0$. Therefore, we have
\begin{eqnarray}
n^{1/2}C\left\{\bTheta_*(\bomega_n)-\bTheta_0\right\}
&=&n^{1/2}C\left\{\frac{\partial\bTheta_*}{\partial\bomega}\Bigm |_{\bomega=\bomega_0}\Delta\bfeta n^{-1/2}+
O(n^{-1})\right\} \nonumber \\
&\rightarrow& \Delta C\left(
\frac{\partial\bTheta_*}{\partial\bomega}\Bigm |_{\bomega=\bomega_0}\right)\bfeta.
\label{secondterm}
\end{eqnarray}
Based on
O'Brien {\em et.\ al.} (2006),
we apply (\ref{secondterm}) in (\ref{truedist}).
Hence we have
the asymptotic distribution of $W$
as noncentral $\chi^2_r(\delta)$,
where the
noncentrality parameter $\delta$ is given by
\begin{eqnarray}\label{non}
\delta&=&\left\{\Delta
C\frac{\partial\bTheta_*}{\partial\bomega}\Bigm
|_{\bomega=\bomega_0} \bfeta\right\}^T
\left\{CI^{-1}_F(\bTheta_0)C^{T}\right\}^{-1}
\left\{\Delta C\frac{\partial\bTheta_*}{\partial\bomega}\Bigm |_{\bomega=\bomega_0}\bfeta\right\},\nonumber\\
&=& \Delta^2 \; \bfeta^T \EE_{\bTheta_0}\{s_G(\bomega_0;Y,\bX)
s^\prime_F(\bTheta_0;Y,\bX)\} I^{-1}_F(\bTheta_0) C^{T}
\left\{CI^{-1}_F(\bTheta_0)C^T\right\}^{-1}
C\nonumber\\
&&
I^{-1}_F(\bTheta_0)
\EE_{\bTheta_0}\{s_F(\bTheta_0;Y,\bX)s^\prime_G(\bomega_0;Y,\bX)\}
\bfeta.
\end{eqnarray}
In the case of a correctly specified model, i.e., ${\cal F} = {\cal G}$,
this reduces to
$$
\delta=\Delta^2
\bfeta^{T}C^{T}
\left\{CI^{-1}_F(\bTheta_0)C^{T}\right\}^{-1}
C \bfeta .
$$
In either case, the asymptotic power of test statistic $W$ is
$$
P\big(\chi^2_r(\delta)>\chi^2_{r,\;\alpha}\big),
$$
where $\chi^2_{r,\;\alpha}$ is the upper $\alpha$
quantile of $\chi^2_r$.

\section{Comparison Between Pairwise Interaction Models
and Diffuse Interaction Models}

\subsection{Diffuse Interaction Model}

Greenland (1983) pointed out that the power of statistical tests
to detect interactions is very low in some commonly encountered
epidemiological situations. We certainly envision low power in
situations where the number of covariates is rather large and only
a very small fraction of all possible (pairwise) interaction terms
really play a role. Gustafson {\em et.\ al.} (2005) proposed
another kind of interaction model, the {\it diffuse interaction}
model, to deal with difficulties caused by a large number of
covariates under pairwise interaction models. The basic form of
this model is best understood in the context of known effect
directions.   For instance, say $Y$ represents a health outcome
(larger values worse), and each $X_{i}$ is a risk factor, scaled
to be nonnegative, such that $E(Y|\bX)$ is known {\em a priori} to
be non-decreasing in each $X_{i}$. Then the DIM form is
\begin{eqnarray}\label{diffuse}
\EE(Y|X_1,\ldots,X_p)&=&
\beta_0+
\left\{\sum_{i=1}^p(\beta_iX_i)^{\lambda}\right\}^{1/\lambda},
\end{eqnarray}
with $\beta_{i} \geq 0$ for $i=1, \ldots, p$.
Note that if $X_{j}=0$ can be interpreted as `absence' of the $j$-th risk factor,
then $\beta_{j}$ can be interpreted as the effect of
$X_{j}$ when all other risk factors are absent, regardless of the value of
$\lambda$.
Assuming normal, homoscedastic errors with $\sigma^{2}=Var(Y| \bX)$,
$(\Beta,\lambda,\sigma^{2})$ comprise the $p+3$ unknown parameters in the
DIM.

To interpret $\lambda$, note first that when $\lambda=1$, (\ref{diffuse}) reduces to the usual additive model.
If $\lambda>1$ though,
then the interaction is antagonistic,
in the sense that for $a<b$,
$\EE(Y|X_j=b,\bX_{(j)}=\bx_{(j)})-\EE(Y|X_j=a,\bX_{(j)}=\bx_{(j)})$
is positive, but
decreasing in each component of $\bx_{(j)}=(x_{1},\ldots, x_{j-1},x_{j+1},\ldots, x_{p})$.
In the special case where each $X_{j}$ is binary to indicate absence or presence of a risk factor, $\lambda>1$ corresponds to
the effect of a particular risk factor diminishing as other risk factors
become present.
Conversely, $\lambda < 1$ corresponds to synergism with the effect growing
as other risk factors become present.
Thus the DIM (\ref{diffuse}) allows for a general tendency for antagonism or synergism
in how multiple risk factors operate on the outcome, without attempting to model fine structure of how
such antagonism or synergism arises.
That is, $\lambda$
controls a one-parameter extension of an additive model which does not single out any particular subset
of risk factors as being less or more responsible for
interaction.

\subsection{Power Comparison between PIM-based and DIM-based Tests}\label{com}

The standard strategy for modelling interactions involves a pairwise-interaction model (PIM).
Again assuming normal, homoscedastic errors,
\begin{eqnarray}\label{standard}
\EE(Y|X_1,\ldots,X_p)&=\beta_0+\sum_{i=1}^p\beta_iX_i+
                      \sum_{i<j}\gamma_{ij}X_iX_j,
\end{eqnarray}
comprises a model with $p(p+1)/2 + 2$ parameters,
$(\Beta,\bgamma,\sigma^{2})$.

To test for departures from additivity then, we could fit the DIM and
test the null that $\lambda=1$, or fit the PIM and test the null that
$\bgamma={\bf 0}$.  In either case, the same null model arises, i.e.,
the additive model with $p+2$ unknown parameters $(\Beta,
\sigma^{2})$.

We can specify a distribution for $X$, values for
$(\Beta , \sigma^{2})$, and a choice of true alternative (DIM or PIM), and then compute the asymptotic power for
both the DIM-based test and the PIM-based test.
When the true alternative is based on DIM, we simply have $
\lambda=1+n^{1/2}\Delta\bfeta$ (since only a single parameter $\lambda$ describes the departure from additivity).
When the true alternative is based on PIM, we must specify the
$p(p-1)/2$ elements of the unit-vector $\bfeta$, i.e., we must specify how the pairwise-interaction coefficients deviate from zero.
Thus investigating the power to detect interactions of PIM form is necessarily more involved than
in the DIM case.

Note that in all cases the quantities needed to determine the
asymptotic power are expectations of squares and
cross products of score vectors for the two models.
Some calculations lack a closed-form due to the particular form of the score vector for the DIM.
The components of this score vector are given in Appendix.
It is the element corresponding to the partial derivative
with respect to $\lambda$ (evaluated at $\lambda=1$) that causes the difficulty
in obtaining an analytical form.
At least in situations where the distribution of $X$ is discrete,
all expectations required
can be calculated via analytic expectation (for $Y|\bX$) and finite summation
for all possible values of $\bX$.
More generally, if $\bX$ follows some continuous distribution, the numerical
integration is required.

\subsection{Detecting Interactions of DIM Form}

\begin{figure}[H]
\centering
\caption{Asymptotic power for a DIM alternative: $X_i,i=1,\dots,9
\stackrel{i.i.d.}{\sim}\mbox{Bernoulli}(0.5)$;
the solid line denotes power of the DIM-based test; the dashed
line denotes power of the PIM-based test.
}
\label{power5}
\includegraphics[scale=0.39,angle=-90]{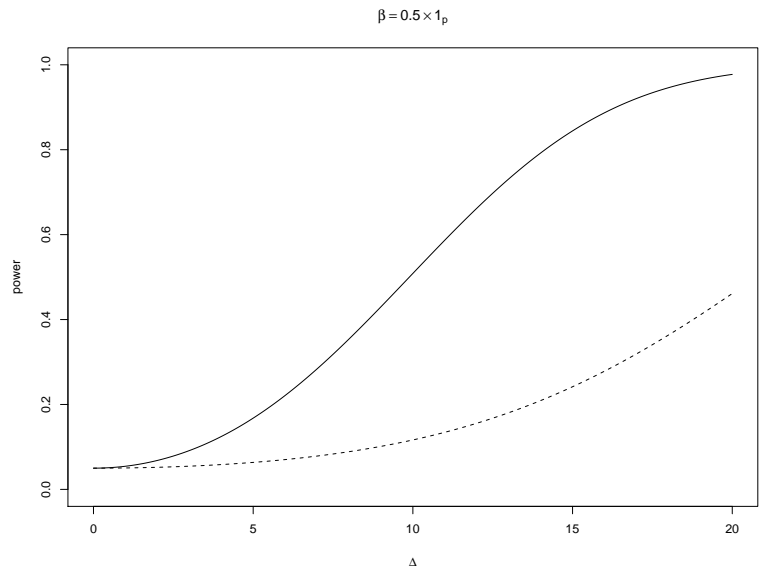}
\end{figure}

As one particular example, say that $p=9$ covariates are
independent and identically distributed as Bernoulli(0.5). Say
that $\Beta^\prime = (0,0.5\times\mathbf{1}^\prime_9)$ and
$\sigma^{2}=1$. Asymptotic power curves (power as a function of
$\Delta$) for the DIM and PIM tests, when the true alternative is
DIM, appear in Figure $\ref{power5}$. As might be anticipated, the
DIM test has substantially higher power, i.e., one does better if
one models the alternative hypothesis correctly.

We find that this conclusion is maintained as we vary the number
and distribution of covariates, and the values of $\Beta$ and
$\sigma$. To some extent we can see this analytically. For
instance, changing $\sigma$ has the same effect of considering a
different value of $\Delta$, as the noncentrality parameter in
($\ref{non}$) is proportional to $\Delta^2/\sigma^2$.

\subsection{Detecting Interactions of PIM Form}

To investigate power when the true alternative follows the PIM, for now we keep the same distribution of $X$
and choice of $(\Beta,\sigma^{2})$ as before,
but must consider different possibilities for $\bfeta$,
the direction in which the pairwise-interaction coefficients deviate from zero.
In an attempt to be somewhat comprehensive, we set up three {\it primary factors} as follows.
Factor 1 is the proportion of entries in $\bfeta$ which are non-zero.
Of the non-zero entries, Factor 2 is the proportion which are positive.
Specializing to the case that all positive entries share the same magnitude and all
negative entries share the same magnitude, Factor 3
is the ratio of the (unique) magnitude of the positive entries to the
(unique) magnitude of the negative entries.
Under this specialization, and given the restriction
$\|\bfeta\|=1$, specifying the three factors does yield a specific value of
$\bfeta$, up to permutation. Here the permutation we choose is
$(0,\ldots,0,1,\ldots,1,-1,\ldots,-1)$.

We consider three levels for each factor, particularly (0.2,0.5,0.8)
for Factor 1 and 2, and (0.5,1,2) for Factor 3.  In the following
plots (Figures 2 through 4), we refer to the three levels of each Factor
as {\it low, medium, high} respectively.

Note first that if the level of Factor 1 is not low and
 the levels of Factors 2 and 3 are both high (or both low), DIM
is be more powerful than PIM. This point is indicated by the superiority of DIM-based
power curves in the two bottom left (right) panels in Figure 2
(Figure 4).
In other words, given a moderate to high presence of pairwise-interaction terms,
if the proportion of the positive (or negative) pairwise terms
with larger magnitudes overwhelms negative (or positive) ones,
i.e., the ``overall" interaction strength leans to synergism (or antagonism), then
DIM works better. However, if the proportions of two opposite directions
are almost equal and the magnitudes of two signs are almost equal as well, DIM
does not work well, as shown in the middle columns of Figures 2 through 4.

Note also that the first column and third column are actually identical in Figure $\ref{clmed}$.
This is caused by the asymptotic power being an even function of $\Delta$, which can be
immediately shown by ($\ref{non}$). Note then the primary factors of (0.2,0.8,1) with
$\Delta>0$ and the primary factors of (0.8,0.2,1) with $\Delta<0$ give the
same value of $\bfeta$.

While Figures 2 through 4 compare the DIM and PIM-based tests
across different settings of the primary factors, there are of
course numerous secondary factors which might be varied as well.
These include the number of covariates and their joint
distribution, as well as the values of the main effect
coefficients.  As one example, Figures 5 through 7 make
comparisons as per Figures 2 through 4, but with the number of
covariates doubled ($p=18$ now). Here we see more of a tendency
for the DIM-based test to compare favourably with the PIM-based
test, as we might expect. R code is available
(www.stat.ubc.ca/$\sim$gustaf) to carry out the power comparison
with primary and secondary factors set to levels desired by the
user.

\section{Discussion}

We view our findings as lending some general support for the
utility of models, such as the DIM, which compromise between the
simplicity of additivity and the flexibility of the PIM.  To
elaborate, we do not claim that the DIM will be highly realistic
across a large range of problems, particularly in the restricted
form considered here (all effect directions known, interactive
behaviour either completely synergistic or completely
antagonistic).  Indeed, Gustafson {\em et.\ al.} (2005) extend the DIM to unknown
directions and Liu (2007) considers more general DIM
forms whereby one group of covariates may have different
interactive behaviour than another.  Even in the simple form
presented here, however, the DIM can capture some coarse structure
of the regression relationship beyond additivity (i.e., a general
tendency for synergistic or antagonistic combination of risk
factors). In contrast, inference in the PIM might be viewed as
attempting to recover fine structure of nonadditive behaviour.
The asymptotic power comparison of PIM and DIM-based tests is
therefore a convenient way to quantify the extent to which coarse
features of non-additivity are more easily detected than fine
features.  It seems interesting that under a true PIM-structure,
enough cohesion in the direction of the pairwise-term coefficients
can render the DIM-based test of non-additivity more powerful than
the PIM-based test.  This matches the applied statistics intuition
that often data will not inform very much about the nature of
nonadditivity, hence a coarse descriptor, such as the single
nonadditivity parameter in the simple DIM, may be
appropriate.

One way in which our stylized treatment of the problem differs from
applied practice is that we have considered the PIM-based test
comparing the additive model with no interactions to the full model
with all possible pairwise interactions.  In practice, particularly
when $p$ is large, one might use a stepwise procedure which
potentially seeks a model with a few pairwise interactions.  Or one
might fit the full model and retain only those pairwise terms with
significant coefficients.  In either case, multiple comparison issues
are at play, and comparison with the DIM-based approach would require
a different strategy than that employed in the present paper.

\section*{References}

\noindent
Le Cam, L. (1960).
Locally asymptotically normal families of distributions.
{\em University of California Publications in Statistics},
{\bf 3}, 37-98.

\noindent
Liu, J. (2007)
Average effects for regression models with misspecifications
and diffuse interaction models, {\em Ph.D. thesis}, Department of Statistics,
University of British Columbia.

\noindent
Greenland, S. (1983).
Tests for interaction in epidemiologic studies:
            a review and a study of power,
{\em Statistics in Medicine},{\bf 2}, 243-251.

\noindent
Gustafson, P. (2001).
On measuring sensitivity to parametric model misspecification.
{\em Journal of the Royal Statistical Society. Series B (Statistical Methodology)},
{\bf 63}, 81-94.

\noindent
Gustafson, P., Kazi, A.M.R. and Levy, A.R. (2005).
Extending logistic regression to model diffuse interactions.
{\em Statistics in Medicine},{\bf 24}, 2089-2104.

\noindent
O'Brien, S.M. and Kupper, L.L. and Dunson, D.B. (2006)
Performance of tests of association in misspecified generalized linear models.
{\em Journal of Statistical Planning and Inference},
{\bf 136}, 3090-3100.

\noindent
White, H. (1982).
Maximum likelihood estimation of misspecified models.
{\em Econometrica},
{\bf 50}, 1-26.

\section*{Appendix}
Let $f_D$ denote the density function of $Y|\bX$ and
$\bTheta=(\Beta^\prime,\lambda,\sigma^2)^\prime$.
Then the score vector for DIM, $s_D(\bTheta;Y,\bX)$, is given by
\begin{eqnarray*}
s_D(\bTheta;Y,\bX)|_{\lambda=1}&=&
\left\{\frac{\partial\log{f_D}}{\partial\bTheta}\right\}\Bigm|_{\lambda=1},\\
&=&\left\{\frac{1}{\sigma^2}(Y-\mu)
\left(\frac{\partial\mu}{\partial\beta_0},\frac{\partial\mu}{\partial\beta_1},
\ldots,\frac{\partial\mu}{\partial\beta_p},\frac{\partial\mu}{\partial\lambda}\right),
\frac{\partial\log f_{D}}{\partial\sigma^2}\right\}\Bigm |_{\lambda=1},
\end{eqnarray*}
where
\begin{eqnarray*}
\mu|_{\lambda=1}&=&\beta_0+\beta_1X_1+\ldots+\beta_pX_p,\nonumber\\
\frac{\partial \mu}{\partial \beta_0}\Bigm |_{\lambda=1}&=&1,\nonumber\\
\frac{\partial \mu}{\partial \beta_j}\Bigm |_{\lambda=1}&=&X_j,j=1,\ldots,p,\nonumber\\
\frac{\partial \mu}{\partial\lambda}\Bigm|_{\lambda=1}&=&
-\left(\sum_{i=1}^p\beta_iX_i\right)\log\left(\sum_{i=1}^p\beta_iX_i\right)
+\sum_{i=1}^p\beta_iX_i\log(\beta_iX_i),\nonumber\\
\frac{\partial\log f_{D}}{\partial\sigma^2}\Bigm|_{\lambda=1}&=&\frac{1}{2\sigma^4}(y-\mu|_{\lambda=1})^2-
\frac{1}{2\sigma^2}.
\end{eqnarray*}

\begin{figure}[H]
\centering \caption{Asymptotic power for PIM alternatives, with
different choices of $\bfeta$ and 9 binary covariates: The three
rows corresponds to three levels of Factor 1, the columns
corresponds to the levels of Factor 2 and Factor 3 is set to be
0.5. Solid lines denote power curves based on diffuse interaction
model fitting and dashed lines denote power curves based on
pairwise interaction model fitting.}\label{cllow}
\includegraphics[scale=0.5,angle=-90]{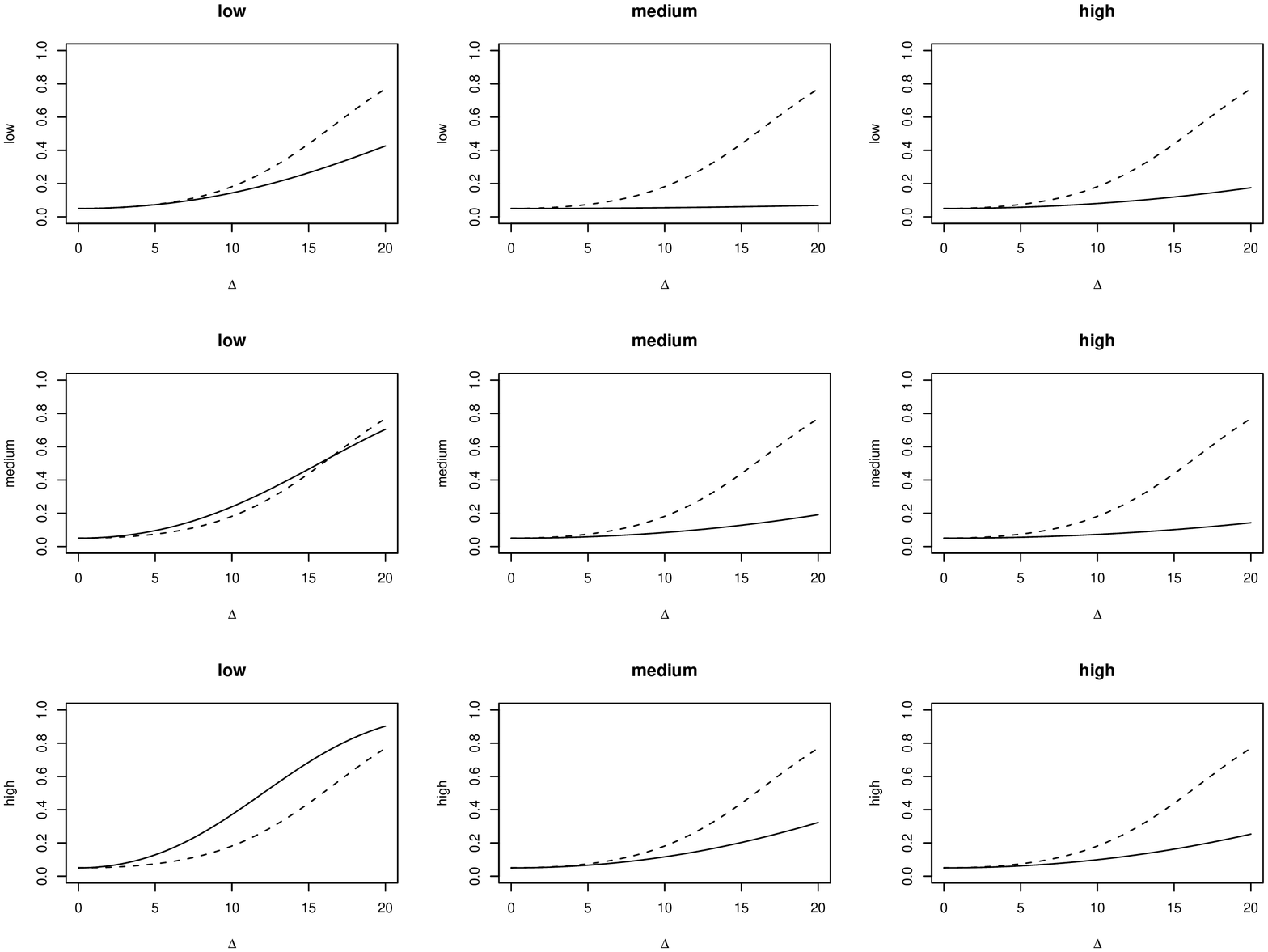}
\end{figure}

\begin{figure}[H]
\centering \caption{Asymptotic power for PIM alternatives, with
different choices of $\bfeta$ and 9 binary covariates: The three
rows corresponds to three levels of Factor 1, the columns
corresponds to the levels of Factor 2 and Factor 3 is set to be 1.
Solid lines denote power curves based on diffuse interaction model
fitting and dashed lines denote power curves based on pairwise
interaction model fitting.}\label{clmed}
\includegraphics[scale=0.5,angle=-90]{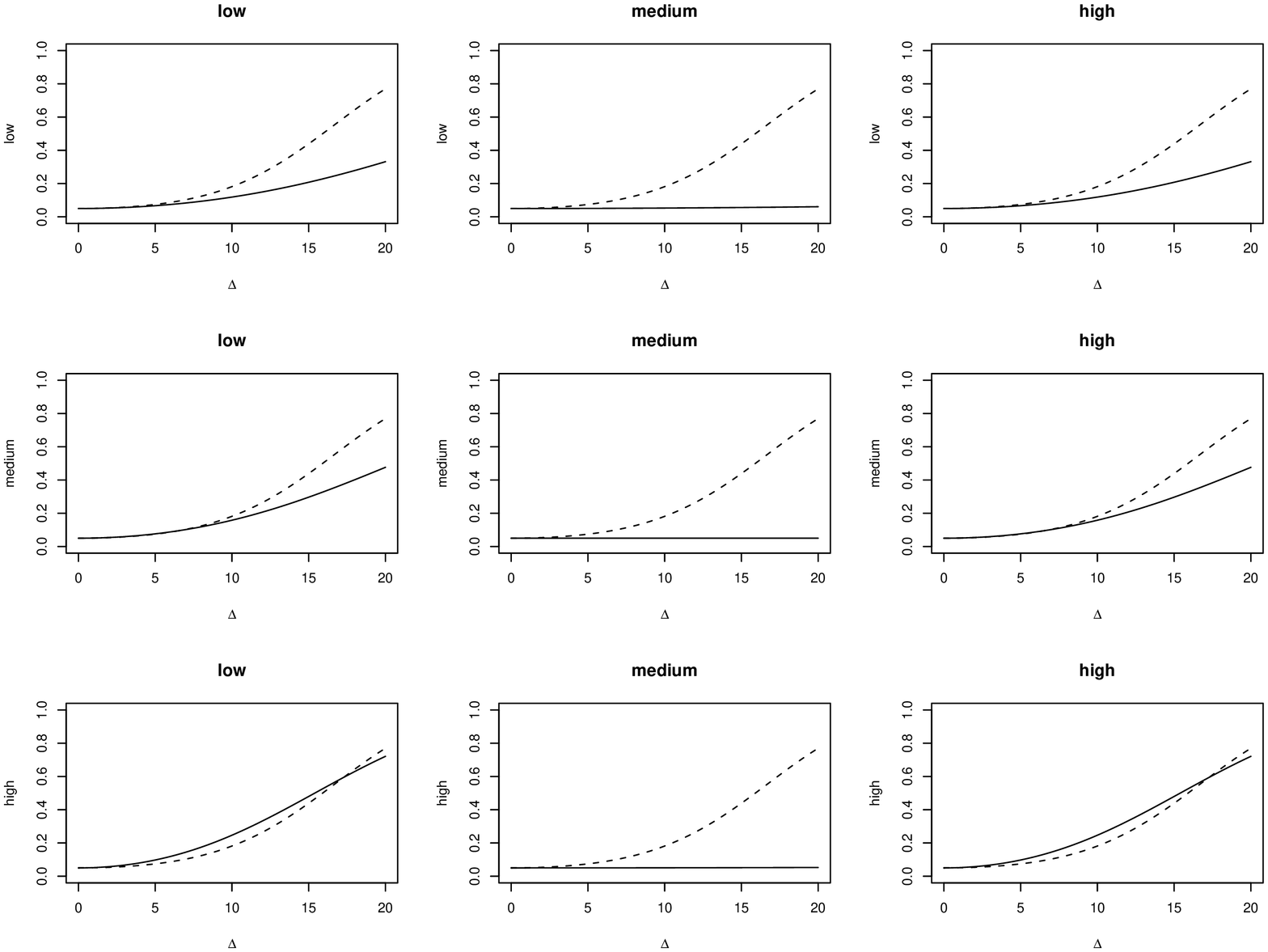}
\end{figure}

\begin{figure}[H]
\centering \caption{Asymptotic power for PIM alternatives, with
different choices of $\bfeta$ and 9 binary covariates: The three
rows corresponds to three levels of Factor 1, the columns
corresponds to the levels of Factor 2 and Factor 3 is set to be 2.
Solid lines denote power curves based on diffuse interaction model
fitting and dashed lines denote power curves based on pairwise
interaction model fitting.}\label{clhigh}
\includegraphics[scale=0.5,angle=-90]{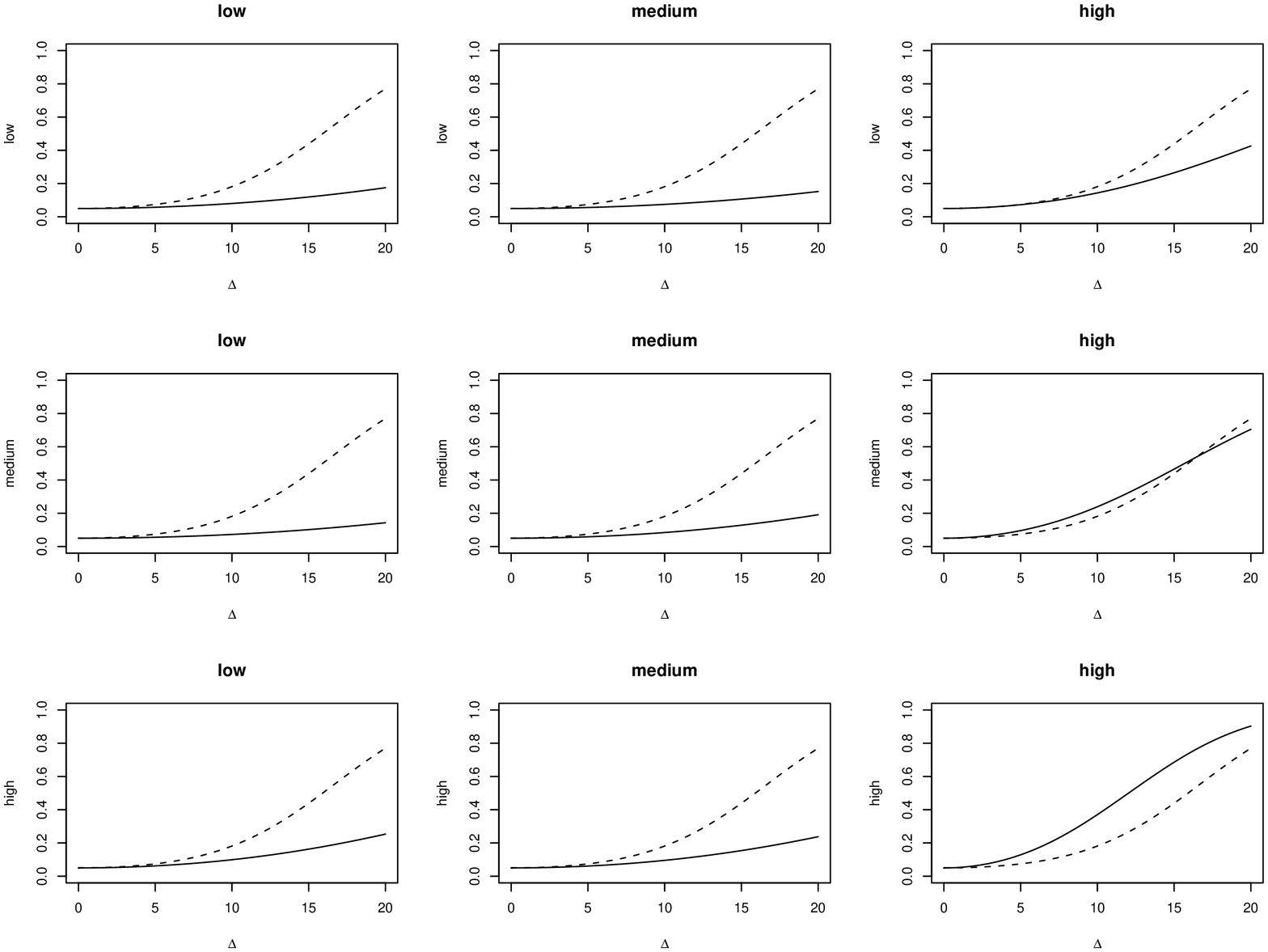}
\end{figure}

\begin{figure}[H]
\centering \caption{Asymptotic power for PIM alternatives, with
different choices of $\bfeta$ and 18 binary covariates: The three
rows corresponds to three levels of Factor 1, the columns
corresponds to the levels of Factor 2 and Factor 3 is set to be
0.5. Solid lines denote power curves based on diffuse interaction
model fitting and dashed lines denote power curves based on
pairwise interaction model fitting.}\label{cllowbp}
\includegraphics[scale=0.5,angle=-90]{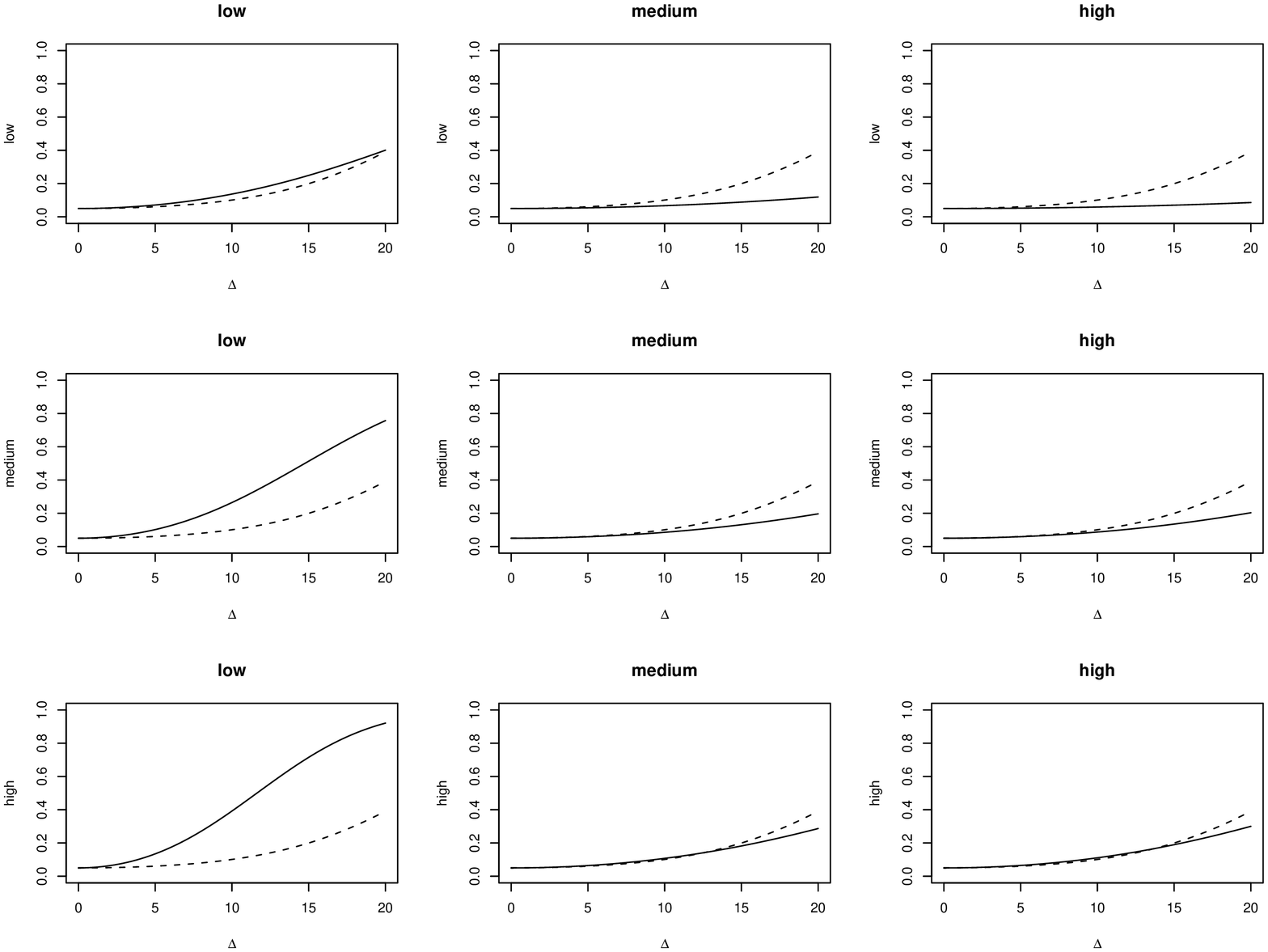}
\end{figure}

\begin{figure}[H]
\centering \caption{Asymptotic power for PIM alternatives, with
different choices of $\bfeta$ and 18 binary covariates: The three
rows corresponds to three levels of Factor 1, the columns
corresponds to the levels of Factor 2 and Factor 3 is set to be 1.
Solid lines denote power curves based on diffuse interaction model
fitting and dashed lines denote power curves based on pairwise
interaction model fitting.}\label{clmedbp}
\includegraphics[scale=0.5,angle=-90]{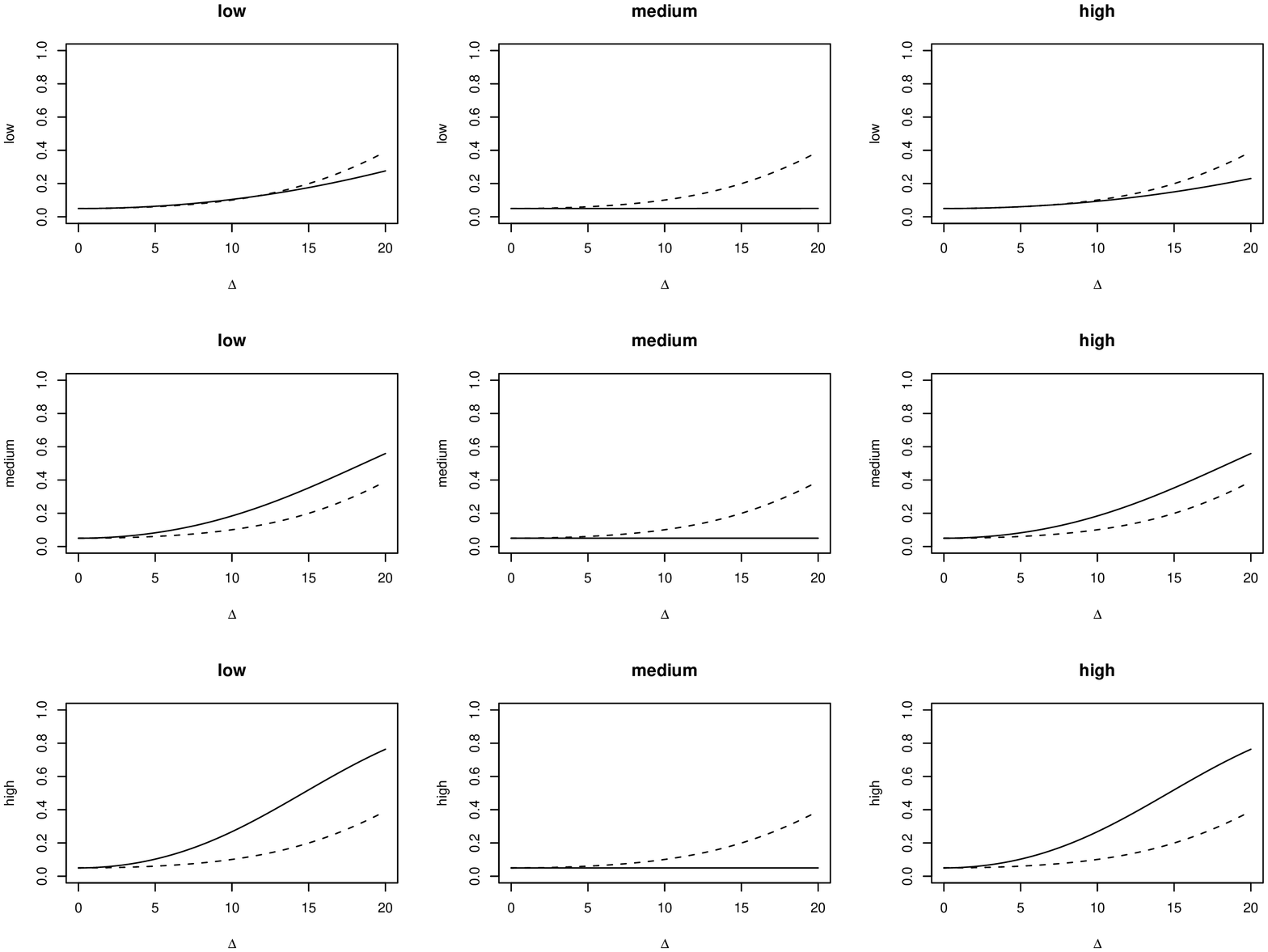}
\end{figure}

\begin{figure}[H]
\centering \caption{Asymptotic power for PIM alternatives, with
different choices of $\bfeta$ and 18 binary covariates: The three
rows corresponds to three levels of Factor 1, the columns
corresponds to the levels of Factor 2 and Factor 3 is set to be 2.
Solid lines denote power curves based on diffuse interaction model
fitting and dashed lines denote power curves based on pairwise
interaction model fitting.}\label{clhighbp}
\includegraphics[scale=0.5,angle=-90]{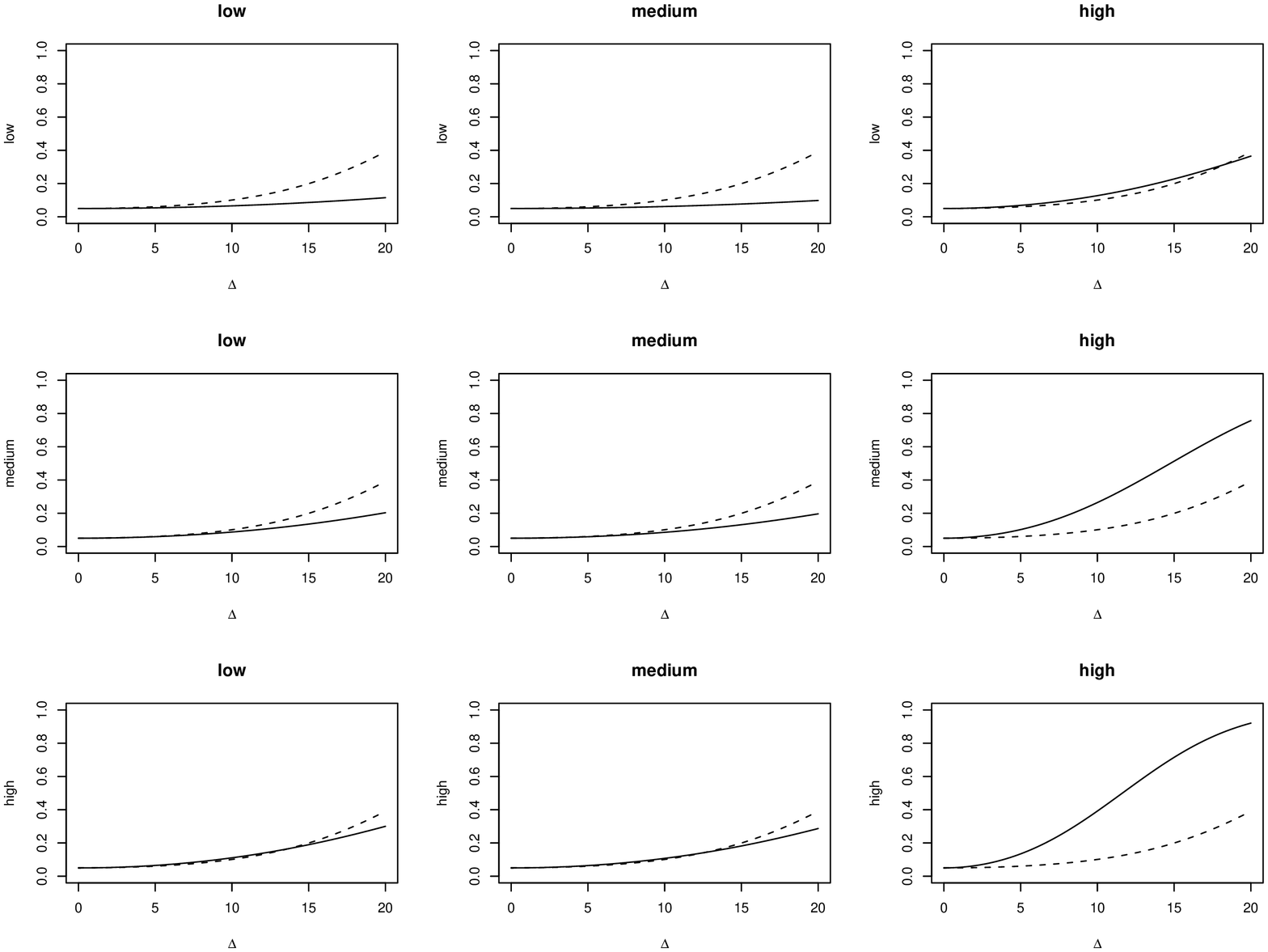}
\end{figure}

\end{document}